\documentclass[12pt,leqno]{article}
\usepackage{amsmath, amsthm, amsfonts, amssymb, color}
\usepackage{mathrsfs}
\usepackage[notcite,notref]{showkeys}
\usepackage[latin1]{inputenc}
\usepackage{bbm}
\usepackage{stmaryrd}
\usepackage{latexsym}
\usepackage{amsfonts}
\usepackage{amsmath,amssymb,amscd}
\usepackage{yfonts}
\usepackage{dsfont}
\usepackage[dvips]{graphicx}
\usepackage{epsfig}
\usepackage{psfrag}

\newtheorem{thm}{Theorem}[section]

\newtheorem{prp}[thm]{Proposition}

\theoremstyle{definition}

\newcommand{\scr}[1]{\mathscr #1}
\definecolor{wco}{rgb}{0.5,0.2,0.3}

\newcommand{\bq}{\begin{eqnarray*}}
\newcommand{\bqn}[1]{\begin{eqnarray}\label{#1}}
\newcommand{\eq}{\end{eqnarray*}}
\newcommand{\eqn}{\end{eqnarray}}

\newcommand{\ttsim}{\raise.17ex\hbox{$\scriptstyle\mathtt{\sim}$}}

\numberwithin{equation}{section} \theoremstyle{remark}

\newcommand{\ua}{\uparrow}

\title{
{\bf  Poincar\'e Inequality for Dirichlet Distributions and Infinite-Dimensional Generalizations}
\footnote{Supported in part by NNSFC(11131003, 11431014), the 985 project, the Laboratory of Mathematical and  Complex Systems,  NSERC, and ANR-STAB-12-BS01-0019.}
}
\author{{\bf Shui Feng$^{b}$, Laurent Miclo$^{d}$,   Feng-Yu Wang$^{a,c}$}\\
\footnotesize{$^a$School of Mathematical Sciences, Beijing Normal University, Beijing 100875, China}\\
 \footnotesize{$^b$Department of Mathematics and Statistics, McMaster University, Hamilton, L8S 4K1, Canada}\\
  \footnotesize{$^c$Department of Mathematics, Swansea University, Singleton Park, SA2 8PP, UK}\\
 \footnotesize{$^d$IMT, UMR 5219, CNRS et Universit\'e Paul Sabatier, 31062 Toulouse cedex 4, France}\\
\footnotesize{  \tttext{wangfy@bnu.edu.cn}; \tttext{laurent.miclo@math.univ-toulouse.fr}; \tttext{shuifeng@univmail.cis.mcmaster.ca}}
}

\begin{document}
\def\tttext#1{{\normalfont\ttfamily#1}}
\def\R{\mathbb R}  \def\ff{\frac} \def\ss{\sqrt} \def\B{\mathbf B}
\def\N{\mathbb N} \def\kk{\kappa} \def\m{{\bf m}}
\def\dd{\delta} \def\DD{\Delta} \def\vv{\varepsilon} \def\rr{\rho}
\def\<{\langle} \def\>{\rangle} \def\GG{\Gamma} \def\gg{\gamma}
  \def\nn{\nabla} \def\pp{\partial} \def\EE{\scr E}
\def\d{\text{\rm{d}}} \def\bb{\beta} \def\aa{\alpha} \def\D{\scr D}
  \def\si{\sigma} \def\ess{\text{\rm{ess}}}
\def\beg{\begin} \def\beq{\begin{equation}}  \def\F{\scr F}
\def\Ric{\text{\rm{Ric}}} \def\Hess{\text{\rm{Hess}}}
\def\e{\text{\rm{e}}} \def\ua{\underline a} \def\OO{\Omega}  \def\oo{\omega}
 \def\tt{\tilde} \def\Ric{\text{\rm{Ric}}}
\def\cut{\text{\rm{cut}}} \def\P{\mathbb P}
\def\C{\scr C}     \def\E{\mathbb E}\def\y{{\bf y}}
\def\Z{\mathbb Z} \def\II{\mathbb I}
  \def\Q{\mathbb Q}  \def\LL{\Lambda}\def\L{\scr L}
  \def\B{\scr B}    \def\ll{\lambda} \def\a{{\bf a}} \def\b{{\bf b}}
\def\vp{\varphi}\def\H{\mathbb H}\def\ee{\mathbf e}\def\x{{\bf x}}
\def\gap{{\rm gap}}\def\PP{\scr P}\def\p{{\mathbf p}}\def\NN{\mathbb N}
\def\cA{\scr A} \def\cQ{\scr Q}

\maketitle
\begin{abstract} For any $N\ge 2$ and $\aa^{(N)}:=(\aa_1,\cdots, \aa_{N+1})\in (0,\infty)^{N+1}$, let $\mu^{(N)}_{\aa}$ be the corresponding Dirichlet distribution on   $\DD:= \big\{ x=(x_i)_{1\le i\le N}\in [0,1]^N:\ \sum_{1\le i\le N} x_i\le 1\big\}.$ We prove  the Poincar\'e inequality
 $$\mu^{(N)}_{\aa}(f^2)\le \ff 1 {\aa_{N+1}} \int_{\DD}\Big\{\Big(1-\sum_{1\le i\le N} x_i\Big) \sum_{n=1}^N   x_n(\pp_n f)^2\Big\}\mu^{(N)}_\aa(\d x)+\mu^{(N)}_{\aa}(f)^2,\   f\in C^1(\DD)$$ and show that the constant $\ff 1 {\aa_{N+1}}$ is sharp.
  Consequently, the associated diffusion process on $\DD$ converges to $\mu^{(N)}_{\aa}$ in $L^2(\mu^{(N)}_{\aa})$ at the exponentially rate $\aa_{N+1}$. The whole spectrum of the generator    is also characterized.  Moreover, the sharp Poincar\'e inequality is extended to the infinite-dimensional setting, and the spectral gap of the corresponding discrete model is derived.     \end{abstract} \noindent

   \vskip 2cm

\section{Introduction}

 Let $N\in \N.$ For any $\aa=(\aa_1,\cdots, \aa_{N+1})\in (0,\infty)^{N+1},$    the Dirichlet distribution $\mu^{(N)}_{\aa}$ with parameter $\aa$ is a probability measure on the set
 $$\DD^{(N)}:= \Big\{ x=(x_i)_{1\le i\le N}\in [0,1]^N:\ \sum_{1\le i\le N} x_i\le 1\Big\}$$ with the density function
 $$\rr(x_1,\cdots, x_N):= \ff{\GG(|\aa|_1)}{\prod_{1\le i\le N+1} \GG(\aa_i)} (1-|x|_1)^{\aa_{N+1}-1}\prod_{1\le i\le N} x_i^{\aa_i-1},\ \ x\in \DD^{(N)},$$ where $|x|_1:= \sum_{1\le i\le N}|x_i|$ for $x\in \R^N.$  Obviously, $\mu^{(N)}_{\aa}$ identifies to the distribution
 $$\tt\mu^{(N+1)}_{\aa}(\d x, \d y) := \mu_{\aa}^{(N)}(\d x) \dd_{1-|x|_1}(\d y)$$ on the space
 $$\nn^{(N+1)}:= \Big\{(x,y)\in [0,1]^{N+1}:\ y+|x|_1=1\Big\}.$$

The Dirichlet distribution arises naturally in Bayesian inference as conjugate priors for categorical distribution and infinite  non-parametric discrete distributions respectively.   They also arise in population genetics describing the distribution of  allelic frequencies (see for instance \cite{ConMoi69,P1,P2}).  In particular, for a population with $N+1$ allelic types, $x_i (1\le i\le N+1)$ stands for the relative frequency of the $i$-th allele among $N+1$ ones.

The Dirichlet distribution possesses many nice properties. We will use the following partition   (or aggregation) property of
$\tt\mu^{(N+1)}_{\aa}$ for $ \aa \in (0,\infty)^{N+1}$. Let   $(X_1, \ldots,X_{N+1})$ have law $\tt\mu^{(N+1)}_{\aa}$,   let $A_1, A_2, \ldots, A_{k+1} $  be a partition of the set $\{1, 2, \ldots, N+1\}$, and set
\[
Y_j=\sum_{r\in A_j}X_r, \  \beta_j=\sum_{r\in A_j}
\alpha_r, \ \ \ j=1,\ldots,k+1.
\]
Then $(Y_1,\ldots, Y_{k+1})$ has law $\tt\mu^{(k+1)}_{\bb}$ with parameters $\bb:=(\beta_1,\ldots,\beta_{k+1})\in (0,\infty)^{k+1}$. We would also like to recall the neutral property of the Dirichlet distribution (\cite{ConMoi69}).  For $(X_1,\cdots, X_N)$ having law $\mu_\aa^{(N)}$,   we define
\[
U_1=X_1,\  U_i =\frac{X_i}{1-X_1-\ldots -X_{i-1}},\ \ \  2\leq i \le N.
\]
Then $U_i$ is a beta random variable with parameters $(\alpha_i, \alpha_{i+1}+\ldots+\alpha_{N+1})$ and $U_1, \ldots, U_N$ are independent.  This leads to the following   representation of the random variable with law $\mu^{(N)}_{\aa}$:
\[
(X_1,X_2,\ldots, X_{N})=\Big(U_1, U_2(1-U_1), \ldots, U_N\prod_{i=1}^{N-1}(1-U_i)\Big).
\]

A well known construction of the Dirichlet distribution is through a P\'olya' urn scheme (\cite{BM73}). More specifically,  consider an urn containing $N+1$ balls of different colors labelled by $1, 2,\ldots, N+1$. The initial mass of the $i$-colored ball is $\alpha_i$. Balls are drawn from the urn sequentially. The chance of a particular colored ball  being selected is proportional to the total mass of that colored  balls inside a urn.   After each selection, the ball is returned with an additional ball of same color and mass one.  The relative weight of different colored balls inside the urn will eventually converge to a Dirichlet vector $(X_1,X_2,\ldots, X_{N+1})$.

To simulate the Dirichlet distribution, several diffusion processes with  this distribution as the stationary distribution  have been proposed and studied. The Wright-Fisher diffusion (see \cite{QQ,Mi1,Mi2,S}) is a diffusion approximation to the Wright-Fisher Markov chain model in population genetics. The evolution mechanism involves mutation and sampling replacement. It is reversible with respect to the Dirichlet distribution.  Exploring the property of right neutrality, a GEM diffusion is introduced in \cite{FW07} and studied further in \cite{FW14}.  This is also a reversible diffusion with Dirichlet distribution as the reversible measure.   A key problem of the study is to estimate the speed for the diffusion process to converge to the  Dirichlet distribution.

The infinite-dimensional generalization of the Dirichlet distribution is Ferguson's Dirichlet process (\cite{Fer73}). It is a random atomic probability measure characterized by the property that its restriction on any finite partition of the state space is a Dirichlet distribution.  The masses of the atoms follow the GEM distribution (see \cite{Ewens04}). The infinite-dimensional generalization of the Wright-Fisher diffusion is the well known Fleming-Viot process with parent independent mutation (\cite{FV79, EK93}) which has an unlabelled version, the infinitely-many-neutral-alleles model  (\cite{EK81}).  The Dirichlet process is the reversible measure of the Fleming-Viot process. Various functional inequalities have been studied to investigate the efficiency of these processes in approximating the Dirichlet processes (\cite{EG93}, \cite{S}).  In particular,  it was shown in \cite{S} that a logarithmic Sobolev inequality holds for the Fleming-Viot process if and only if the state space has finite dimension.  The Wasserstein diffusion studied in \cite{RS09} and \cite{DS09} is closely related to the Dirichlet process with state space $[0,1]$.   However,  the exact convergence rate is not yet known for these processes to approximate the Dirichlet distribution. 

In this paper, we will construct and study a new class of processes with Dirichlet distribution and Dirichlet processes as reversible measures. In comparison with the Wright-Fisher Markov chain model where sampling replacement occurs between any pair of individuals, our finite dimensional Markov chain model only allows sampling replacement between individuals in one group and individuals in another group. This reduced sampling scheme bridges the gap between the independent systems and the pairwise sampling models. A complete understanding of these models will provide a whole picture of the roles played by different evolutionary forces.  

The main contributions  of this paper are the explicit identification of the whole spectrum of the finite-dimensional diffusions, the construction of an infinite-dimensional diffusion process with Dirichlet process on a countable state space as the reversible measure, the establishment of sharp Poincar\'e inequalities for both the finite and infinite dimensional diffusions, and the construction of Markov chain models.  In particular, we found the exact exponential  rate of  the underlying diffusion processes and Markov chains to converge to the Dirichlet distribution  and its 
infinite-dimensional generalization.  

An outline of development of the paper is as follows. In Section 2, we collect the main results. The sharp Poincar\'e inequality for the finite-dimensional diffusion is proved in Section 3. The key step in the proof is to link the eigenvalues of the diffusion generator with a finite matrix.  The whole spectrum is obtained in Section 4.  Section 5 contains the construction of the infinite-dimensional diffusion and the establishment of the corresponding  sharp Poincar\'e inequality. Finally, in Section 6  we introduce the discrete Markov chain model involving immigration, emigration and sampling, which approximates the diffusion process solving \eqref{E1}.  This not only provides the genetic link to the diffusion model  but also opens the door for the study of models with other sampling mechanisms.

\section{Main Results}

The diffusion process studied in this paper first appeared in \cite[(2.44)]{Jac01} (see also \cite{BR}) and  solves the following SDE on $\DD^{(N)}$:
\beq\label{E1} \d X_i(t)= \big\{\aa_i (1-|X(t)|)- \aa_{N+1}X_i(t)\big\}\d t + \ss{2 (1-|X(t)|_1)X_i(t)}\,\d B_i(t),\  \ 1\le i\le N,\end{equation}
where  $B(t):=(B_1(t),\cdots, B_N(t))$ is the $d$-dimensional Brownian motion.

We will show that  the Markov semigroup $P_t^\aa$ associated to \eqref{E1} is symmetric in $L^2(\mu^{(N)}_{\aa})$; that is,
 \beq\label{E2} \int_{\DD^{(N)}} fL^{(N)}_{\aa}g \d\mu^{(N)}_{\aa}= \int_\DD^{(N)} gL_\aa^{(N)} f \d\mu^{(N)}_{\aa},\ \ f,g\in C^2(\R^N) \end{equation} holds for
 $$L^{(N)}_\aa(x):= \sum_{1\le n\le N} \Big(x_n (1-|x|_1)  \pp_n^2 + \big\{\aa_n (1-|x|_{1})-\aa_{N+1}x_n\big\}\pp_n\Big) $$
being the generator of $P_t^\aa$, where   and $\pp_n:=\ff{\pp}{\pp x_n}$.
So, $(L_\aa, C^2({\DD^{(N)}}))$ is closable in $L^2(\mu^{(N)}_{\aa})$ and its closure $(L_\aa,\D(L_\aa))$ is a negative definite self-adjoint operator.

Moreover, since
$$L_{\aa}^{(N)}(fg)(x)= (fL_\aa^{(N)}g+ gL_{\aa}^{(N)} f)(x) + 2  (1-|x|_1) \sum_{n=1}^N   x_n\{(\pp_n f)(\pp_ng)\}(x),$$
\eqref{E2} implies the integration by parts formula
\beq\label{E2'} \beg{split} -\int_{\DD^{(N)}} fL_{\aa}^{(N)} g \d\mu_\aa^{(N)}&=  \int_{\DD^{(N)}}\Big\{ (1-|x|_1) \sum_{n=1}^N   x_n\{(\pp_n f)(\pp_ng)\}(x)\Big\} \mu_\aa^{(N)}(\d x)\\
&=:\EE_\aa^{(N)}(f,g),\ \ \ f,g\in C^2(\DD^{(N)}).\end{split}\end{equation}
Therefore, $(\EE_\aa^{(N)}, C^2(\DD^{(N)}))$ is closable in $L^2(\mu_\aa^{(N)})$ whose closure
  $(\EE_\aa^{(N)}, \D(\EE_\aa^{(N)}))$ is a symmetric Dirichlet form on $L^2(\mu_\aa^{(N)})$, and it is easy to see  that  this Dirichlet form  is associated to  the Markov semigroup $P_t^\aa$.

Finally,   the spectral gap of $L_\aa^{(N)}$ is characterized as
$$\gap(L_\aa^{(N)})= \inf\Big\{\EE_\aa^{(N)}(f,f):\ f\in \D(\EE_\aa^{(N)}), \mu^{(N)}_{\aa}(f)=0, \mu^{(N)}_{\aa}(f^2)=1\Big\}.$$
It is known that when $N=1$ we have $\gap(L_\aa^{(N)})= \aa_1+\aa_2$, see e.g. \cite{S}. So, in the following result   we only consider $N\ge 2$.

\beg{thm}\label{T1.1} Let $N\ge 2$. Then $P_t^\aa$ is symmetric in $L^2(\mu_\aa^{(N)})$ and its generator has spectral gap
$\gap(L_\aa^{(N)})= \aa_{N+1}.$ Consequently, $P_t^\aa$ converges to $\mu^{(N)}_{\aa}$ exponentially fast in $L^2(\mu^{(N)}_{\aa}):$  $$\|P_t^\aa-\mu^{(N)}_{\aa}\|_{L^2(\mu^{(N)}_{\aa})}\le \e^{-\aa_{N+1}t},\ \ t\ge 0,$$ and the sharp Poincar\'e inequality for $(\EE_\aa^{(N)},\D(\EE_\aa^{(N)}))$ is
$$\mu^{(N)}_{\aa}(f^2)\le \ff 1 {\aa_{N+1}}\, \EE_\aa^{(N)}(f,f),\ \ f\in \D(\EE_\aa^{(N)}), \mu_\aa^{(N)}(f)=0.$$\end{thm}

Next, we extend this result to the infinite-dimensional setting. Consider the infinite-dimensional simplex
$$\DD^{(\infty)} := \Big\{x\in [0,1]^\N:\ |x|_1 =\sum_{i=1}^\infty  x_i \le 1\Big\},$$ which is equipped with the $L^1$-metric $  |x-y|_1.$ Let $\aa\in (0,\infty)^\N$ with $|\aa|_1 =\sum_{i=1}^\infty\aa_i<\infty$, and let $\aa_\infty>0$ which refers to $\aa_{N+1}$ in the finite-dimensional case as $N\to\infty$. Let  $$\aa^{(n)}=\Big(\aa_1,\cdots, \aa_{n-1},\sum_{i\ge n}\aa_i, \aa_\infty\Big)\in (0,\infty)^{n+1},\ \ \ n\ge 1.$$ Then for any $n\ge 1$,
 $$\mu_{\aa,\aa_\infty}^{(n)}(\d x):= \mu_{\aa^{(n)}}^{(n)}(\d x_1,\cdots, \d x_n) \prod_{i=n+1}^\infty \dd_0(\d x_i)$$ is a probability measure
 on $\DD^{(\infty)}.$ We will prove that when $n\to \infty$ these measures converges weakly to a probability measure $\mu_{\aa,\aa_\infty}^{(\infty)}$ on $\DD^{(\infty)}$, which is the infinite-dimensional generalization of Dirichlet distribution with parameters $(\aa,\aa_\infty)$. The following result extends Theorem \ref{T1.1} to the case for $N=\infty$.

\beg{thm}\label{T1.2} Let $\aa\in (0,\infty)^\N$ with $|\aa|_1 <\infty$ and let $\aa_\infty>0$.
\beg{enumerate}
\item[$(1)$] The sequence $\{\mu^{(n)}_{\aa,\aa_\infty}\}_{n\ge 1}$ converges weakly to a probability measure $\mu_{\aa,\aa_\infty}^{(\infty)}$ on $\DD^{(\infty)}$.
\item[$(2)$] The form
$$\EE^{(\infty)}_{\aa,\aa_\infty}(f,g):=\int_{\DD^{(\infty)}}\Big\{(1-|x|_1) \sum_{n=1}^\infty   x_n(\pp_n f)\pp_n g\Big\}(x) \mu_{\aa,\aa_\infty}^{(\infty)}(\d x),\ \ f,g\in \scr FC^1$$ is closable in $L^2(\mu_{\aa,\aa_\infty}^{(\infty)})$ whose closure is a symmetric Dirichlet form. The    generator $(L_{\aa,\aa_\infty}^{(\infty)},\D(L_{\aa,\aa_\infty}^{(\infty)}))$ of the Dirichlet form satisfies $\scr FC^2\subset \D(L_{\aa,\aa_\infty}^{(\infty)})$ and
$$L_{\aa,\aa_\infty}^{(\infty)} f(x)=\sum_{n=1}^\infty \Big(x_n (1-|x|_1)  \pp_n^2f(x) + \big\{\aa_n (1-|x|_{1})-\aa_{\infty}x_n\big\}\pp_nf(x)\Big),
\ \ f\in \scr FC^2.$$
\item[$(3)$] The generator $L_{\aa,\aa_\infty}^{(\infty)}$ has spectral gap $\gap(L_{\aa,\aa_\infty}^{(\infty)})= \aa_\infty.$ Consequently, the associated Markov semigroup $P_t^{\aa,\aa_\infty}$ converges to $\mu_{\aa,\aa_\infty}^{(\infty)}$ exponentially fast in $L^2(\mu^{(\infty)}_{\aa,\aa_\infty}):$  $$\|P_t^{\aa,\aa_\infty}-\mu_{\aa,\aa_\infty}^{(\infty)}\|_{L^2(\mu_{\aa,\aa_\infty}^{(\infty)})}\le \e^{-\aa_{\infty}t},\ \ t\ge 0,$$ and the sharp Poincar\'e inequality  is
$$\mu_{\aa,\aa_\infty}^{(\infty)}(f^2)\le \ff 1 {\aa_{\infty}}\, \EE_{\aa,\aa_\infty}^{(\infty)}(f,f),\ \ f\in \scr FC^1,\mu_{\aa,\aa_\infty}^{(\infty)}(f)=0.$$
\end{enumerate} \end{thm}

Finally, the next result shows that the diffusion process generated by $L_{\aa,\aa_\infty}^{(\infty)}$ is  the weak limit of the $L_{\aa,\aa_\infty}^{(n)}$-diffusion process as  $n\to\infty$, where
$$L_{\aa,\aa_\infty}^{(n)}:= \sum_{i=1}^n \Big\{\Big[\aa_i\Big(1-\sum_{i=1}^n x_i\Big)-\aa_\infty x_i\Big]\pp_i + 2 \Big(1-\sum_{i=1}^n x_i\Big)x_i \pp_i^2\Big\}.$$  For any $x\in \DD^{(\infty)}$ and $T>0$, let $P_{x,T}^{(n)}$ be the distribution of the diffusion process generated by $L_{\aa,\aa_\infty}^{(n)}$ with initial point $x^{(n)}:=\big(x_1,\cdots, x_{n-1}, \sum_{j\ge n} x_j\big)$. Embedding $\DD^{(n)}$ into $\DD^{(\infty)}$ by setting $z_i=0$ for $z\in \DD^{(n)}$ and $i\ge n+1$, we regard $P_{x,T}^{(n)}$ as a probability measure on $\OO_T:=C([0,T];\DD_\infty)$ equipped with the uniform norm $\|\xi\|_{1,\infty}:= \sup_{t\in [0,T]}|\xi(t)|_1.$

\beg{thm}\label{T1.3} For any $x\in \DD^{(\infty)}$ and $T>0$, $P_{x,T}^{(n)}$ converges weakly to a probability measure $P_{x,T}^{(\infty)}$ on
$\OO_T$. Moreover,  $P_{x,T}^{(\infty)}$ solves the martingale problem of $L_{\aa,\aa_\infty}^{(\infty)}$: for any $f\in \F C^2$,   the coordinate process $
X(t)(\oo):=\oo(t)$ and the natural filtration $\F_t:=\si(\oo_s:\ s\in [0, t]), $
$$f(X(t))-\int_0^t L_{\aa,\aa_\infty}^{(\infty)} f(X(s))\d s,\ \ t\in [0,T]$$ is a  martingale under $P_{x,T}^{(\infty)}.$\end{thm}

\section{Proof of Theorem \ref{T1.1}}
We first prove \eqref{E2} which implies the symmetry of $P_t^\aa$ in $L^2(\mu_\aa^{(N)}).$
Since smooth functions on $\DD^{(N)}$ are uniformly approximated by polynomials up to second order derivatives,
it suffices to consider
$f,g\in  \scr P_\infty$, the set of all polynomials on $\DD^{(N)}$. Let
 $$A_\aa^{(n)}= x_n (1-|x|_{1}) \pp_n^2 + \big\{\aa_n (1-|x|_{1})-\aa_{N+1}x_n\big\}\pp_n,\ \ 1\le n\le N.$$  Then \eqref{E2} follows from
 \beq\label{E3} \int_\DD^{(N)} \Big(\prod_{1\le i\le N} x_i^{p_i}\Big) A_\aa^{(n)} \Big(\prod_{1\le i\le N} x_i^{q_i}\Big)  \mu^{(N)}_{\aa}(\d x)
 = \int_\DD^{(N)} \Big(\prod_{1\le i\le N} x_i^{q_i}\Big) A_\aa^{(n)} \Big(\prod_{1\le i\le N} x_i^{p_i}\Big) \mu^{(N)}_{\aa}(\d x) \end{equation} for $p_i,q_i\in \Z_+, 1\le i\le N.$ Letting $p_{N+1}=q_{N+1}=0$ and   $C=  \ff{\GG(|\aa|_1)}{\prod_{1\le i\le N+1} \GG(\aa_i)},$ and simply denote
 $x_{N+1}=1-|x|_1$, we have
 \beg{equation*}\beg{split} &\int_{\DD^{(N)}} \Big(\prod_{1\le i\le N} x_i^{p_i}\Big) A_\aa^{(n)} \Big(\prod_{1\le i\le N} x_i^{q_i}\Big)
   \mu^{(N)}_{\aa}(\d x)\\
  &=  C\int_{\DD^{(N)}} \Big(\prod_{1\le i\ne n \le N+1} x_i^{p_i+q_i+\aa_i-1}\Big)x_n^{p_n+\alpha_n-1}A_\aa^{(n)} x_n^{q_n} \d x\\
 &=  Cq_n \bigg\{(q_n+\aa_n-1) \int_{\DD^{(N)}} \Big(\prod_{1\le i\ne n \le N+1} x_i^{p_i+q_i+\aa_i-1}\Big)
x_{N+1}x_n^{p_n+q_n+\aa_n-2}\d x \\
&\qquad \qquad\qquad\qquad\qquad -\aa_{N+1} \int_{\DD^{(N)}} \Big(\prod_{1\le i  \le N+1} x_i^{p_i+q_i+\aa_i-1}\Big)\d x\bigg\}\\
&= \ff{C q_n\prod_{1\le i\ne n\le N+1} \GG(\aa_i+p_i+q_i)}{ \GG(\sum_{1\le i\le N+1} (\aa_i+p_i+q_i))}\\
&\quad\times \Big((q_n+\aa_n-1)\GG(\aa_{N+1}+1)\GG(p_n+q_n+\aa_n-1) -\aa_{N+1}\GG(\aa_{N+1})\GG(p_n+q_n+\aa_n)\Big)\\
&= - \ff{C \GG(\aa_{N+1}+1)\prod_{1\le i\ne n\le N+1} \GG(\aa_i+p_i+q_i)}{ \GG(\sum_{1\le i\le N+1} (\aa_i+p_i+q_i))} p_nq_n\GG(p_n+q_n+\aa_n-1),\end{split}\end{equation*} where the last step is due to the identity $\GG(s+1)=s\GG(s), s>0.$ Since
the result is symmetric in $(p_n, q_n)$, it implies \eqref{E3}.

\

For any $d\in \N$, let $\PP_{d}$ be the space of all polynomials  in $\scr P_\infty$ whose total degrees are less than or equal to $d$. Let
$\PP_{0, d}=\{f\in \PP_d: \mu^{(N)}_{\aa}(f)=0\}.$ It is well known that $\PP_{\infty}:= \cup_{d\ge 1} \PP_{d}$ is dense in $C_b^1({\DD^{(N)}})$, so that
$\PP_{0,\infty}:=\cup_{d\ge 1} \PP_{0,d}$ is dense in $$\D_0:=\{f\in \D(\EE_\aa^{(N)}):\mu^{(N)}_{\aa}(f)=0\}$$ under the Sobolev norm $\|f\|_{1,2}:= \ss{\mu_\aa^{(N)}(f^2) + \EE_{\aa}^{(N)}(f,f)}$\ .

To characterize $\gap(L_\aa^{(N)})$, we make the spectral decomposition  of $L_\aa^{(N)}$ in terms of the degree of polynomials. Obviously, every $\PP_{0,d}$ is an invariant space of $L_\aa^{(N)}$. Let $\scr Q_1= \PP_{0,1}$ and
$$\scr Q_d=\big\{f\in \PP_{0,d}: \mu^{(N)}_{\aa}(fg)=0 \text{\ for\ all}\  g\in \PP_{d-1}\big\},\ \ d\ge 2.$$ Then, by the symmetry of $L_\aa^{(N)}$ in $L^2(\mu^{(N)}_{\aa})$, every $\scr Q_d$ is an invariant space of $L_\aa^{(N)}$ as well. Thus,   letting $\pi_d: \PP_{\infty} \to \PP_{d}$ be the orthogonal projection with respect to the inner product in $L^2(\mu^{(N)}_{\aa})$, we have
\beq\label{P1} L_\aa^{(N)} \pi_{d}f = \pi_{d} L_\aa^{(N)} f,\ \ d\ge 1, f\in \PP_\infty.\end{equation} Therefore,  to characterize the spectrum of $L_\aa^{(N)}$ it suffices to consider that of    $L_\aa^{(N)}|_{\scr Q_i},$ the restriction  of $L_\aa^{(N)}$ on  $\scr Q_i,$ for every $ i\ge 1$.

Let $d\ge 2$.  To characterize the spectrum of $L_\aa^{(N)}|_{\scr Q_d},$ let $$K_d=\Big\{k=(k_1,\cdots, k_N)\in \Z_+^N:\ \sum_{1\le i\le N} k(i)=d\Big\}.$$ For any $k\in K_d$, let
$x^k= \prod_{1\le i\le N} x_i^{k_i}.$ Then
\beq\label{Q} \scr Q_d= \bigg\{\sum_{k\in K_d} c_k x^k- \pi_{d-1} \sum_{k\in K_d} c_k x^k:\  c:=(c_k)_{k\in K_d}\in \R^{K_d}
\bigg\}.\end{equation}
We define the $K_d\times K_d$-matrix $M_d$ by letting
$$M_d(k,k')=\begin{cases} d\aa_{N+1} +\sum_{1\le n\le N} (k_n+\aa_n-1)k_n, & \text{if}\ k=k',\\
(k_n+\aa_n)(k_n+1), & \text{if}\ k'= k+e_n-e_m, 1\le n\ne m\le N,\\
0, & \text{otherwise},\end{cases}$$ where $\{e_n\}_{1\le n\le N}$ is the canonical orthonormal  basis on $\R^N$.  We first identify eigenvalues of $L_\aa^{(N)}|_{\scr Q_d}$ with those of $M_d$.

\beg{lem}\label{L1} For any $d\ge 2$, $\ll$ is an eigenvalue of $-L_\aa^{(N)}|_{\scr Q_d}$ if and only if it is an eigenvalue of $M_d$. Consequently,
$-L_\aa^{(N)}|_{\scr Q_d}\ge (d \aa_{N+1}) I_{\scr Q_d}$, where $I_{\scr Q_d}$ is the identity operator on $\scr Q_d.$ \end{lem}

\beg{proof} (1) Let $\ll$ be an eigenvalue of $-L_\aa^{(N)}$ on $\scr Q_d$. By \eqref{Q} and \eqref{P1}, there exists $0\ne c\in \R^{K_d}$ such that
\beq\label{W1}  \sum_{k\in K_d} c_k (L_\aa^{(N)} x^k  - \pi_{d-1} L_\aa^{(N)} x^k)=-\ll \sum_{k\in K_d} c_k ( x^k  - \pi_{d-1}  x^k).\end{equation}
Obviously,
\beg{equation*}\beg{split}  & L_\aa^{(N)} x^k -\sum_{1\le n\le N} (x_n \pp_n^2+\aa_n \pp_n) x^k \\
&= -  \bigg( \sum_{1\le n,m \le N}x_nx_m\pp_n^2  x^k
+\sum_{1\le n,m\le N}x_m\alpha_n\pp_n x^k +\alpha_{N+1}\sum_{1\le n\le N}x_n\pp_n x^k \bigg) \\
%&= -\bigg(\sum_{1\le n,m\le N}x_nx_mk_n(k_n-1) x^{k-2e_n}
%+\sum_{1\le n,m\le N}x_m\alpha_n k_nx^{k-e_n}+\alpha_{N+1}\sum_{1\le n\le N}x_nk_nx^{k-e_n}\bigg)\\
&= -\bigg(\sum_{n,m\le N}k_n(k_n-1) x^{k-e_n+e_m}
+\sum_{1\le n,m\le N}\alpha_n k_nx^{k-e_n+e_m}+\alpha_{N+1}\sum_{1\le n\le N}k_nx^{k}\bigg)\\
&= -\bigg(\sum_{1\le n,m\le N}k_n(k_n-1) x^{k-e_n+e_m}
+\sum_{1\le n,m\le N}\alpha_n k_nx^{k-e_n+e_m}+d\alpha_{N+1}x^{k}\bigg).\end{split}\end{equation*}
By the change of variables $k':= k-e_n+e_m$,  we obtain
\beg{equation*}\beg{split}  &\sum_{k\in K_d} c_k \sum_{1\le n,m\le N}\alpha_n k_nx^{k-e_n+e_m}\\
&= \sum_{k\in K_d} c_k \sum_{1\le n\ne m\le N}\alpha_n k_nx^{k-e_n+e_m}+\sum_{k\in K_d} c_k \sum_{1\le n\le N}\alpha_n k_nx^{k}\\
&=\sum_{k\in K_d}  \sum_{1\le n\ne m\le N}c_{k'+e_n-e_m}\alpha_n (k'+e_n-e_m)(n)x^{k'}+\sum_{k\in K_d} c_k \sum_{1\le n\le N}
\alpha_n k_nx^{k}\\
&= \sum_{k\in K_d} \sum_{1\le n\ne m\le N}\alpha_n (k_n+1)c_{k+e_n-e_m}x^{k}+\sum_{k\in K_d}  \sum_{1\le n\le N}\alpha_n k_nc_k x^{k}.\end{split}\end{equation*} Similarly,
\beg{equation*}\beg{split}  &\sum_{k\in K_d} c_k \sum_{1\le n,m\le N}  k_n(k_n-1)x^{k-e_n+e_m}\\
&= \sum_{k\in K_d} \sum_{1\le n\ne m\le N} k_n(k_n+1) c_{k+e_n-e_m} x^k + \sum_{k\in K_d} c_k \sum_{1\le n\le N} k_n(k_n-1)x^k.\end{split}\end{equation*} Combining these together leads to
\beq\label{W2}   \sum_{k\in K_d} c_k  L_\aa^{(N)} x^k = \sum_{k\in K_d} c_k \sum_{1\le n\le N} (x_n \pp_n^2+\dd_n\pp_n) x^k- \sum_{k,k'\in K_d} M_d(k,k')c_{k'} x^k.\end{equation}
Substituting this into \eqref{W1}, we arrive at
$$\sum_{k\in K_d} (M_d c)_k x^k=\ll \sum_{k\in K_d} c_kx^k +\p_{d-1}(x)$$ for some $\p_{d-1}\in \PP_{d-1}.$
Therefore, $M_d c= \ll c$, i.e. $\ll$ is an eigenvalue of $M_d$.

(2) On the other hand, if $\ll$ is an eigenvalue of $M_d$, then there exists $  c\in \R^{K_d}\setminus \{0\}$ such that $M_d c=\ll c$. Let
$$f(x)= \sum_{k\in K_d} c_k x^k - \pi_{d-1} \sum_{k\in K_d} c_k x^k.$$ It follows from $M_d c=\ll c$ and \eqref{W2} that
$$L_\aa^{(N)} f= \tt \p_{d-1}-\ll f $$ holds for some $\tt \p_{d-1}\in \PP_{d-1}.$ Since $f\in \scr Q_d$ which is orthogonal to $\PP_{d-1}$, this and \eqref{P1} implies
$$L_\aa^{(N)} f= (1-\pi_{d-1}) L_\aa^{(N)} f = -\ll (1-\pi_{d-1}) f=-\ll f.$$ So, $\ll$ is an eigenvalue of $L_\aa^{(N)}$ on $\scr Q_d$.

(3) Finally, since eigenvalues of $-L_\aa^{(N)}$ are nonnegative,  (2) implies
 that eigenvalues of $\tt M_d:=M_d-d\aa_{N+1}I_{K_d\times K_d}$ is larger than or equal to $-d\aa_{N+1}.$ On the other hand, from the definition of $M_d$ we   see that
$\tt M_d$ does not depend on $\aa_{N+1}$. So, letting $\aa_{N+1}\downarrow 0$ and noting that $M_d\ge 0$,  we conclude that eigenvalues of $\tt M_d$ are non-negative. Therefore,
eigenvalues of $M_d$ are larger than or equal to $d \aa_{N+1}.$ Combining this with (1)   we obtain  $-L_\aa^{(N)}|_{\scr Q_d}\ge (d \aa_{N+1}) I_{\scr Q_d}.$
\end{proof}

\beg{proof}[Proof of Theorem \ref{T1.1}] By Lemma \ref{L1}, it suffices to prove that the smallest eigenvalue of $-L_\aa^{(N)}|_{\scr Q_1}$ is $\aa_{N+1}.$
To this end, we take $\theta_i=(\theta_{ij})_{1\le j\le N}\in \R^N ( 1\le i\le N-1)$ such that
$$\sum_{k=1}^N \theta_{ik}\aa_k=0,\ \ \sum_{k=1}^N \theta_{ik}\theta_{jk}\aa_k=\dd_{ij},\ \ 1\le i,j\le N-1.$$
So, $\{\theta_i\}_{i=1}^N$ is a basis of $\R^{N-1}$.
Let
\beg{equation*}\beg{split} &u_i(x)= \sum_{j=1}^N \theta_{ij} x_j,\ \ 1\le i\le N-1;\\
& u_N(x)= \sum_{k=1}^N x_k -\ff{\tt\aa}{|\aa|_1},\ \ \tt\aa:=|\aa|_1-\aa_{N+1}=\sum_{k=1}^N\aa_k.\end{split}\end{equation*}
We intend to  prove that $\{u_i\}_{1\le i\le N}$ is an orthogonal basis of $\scr Q_1$ with respect to the inner product $\<f,g\>_\aa^{(N)}:= \mu^{(N)}_{\aa}(fg)=\int_{\DD^{(N)}} fg\d\mu^{(N)}_{\aa},$ and $L_\aa^{(N)} u_N=-|\aa|_1 u_N$ while $L_\aa^{(N)} u_i=-\aa_{N+1}u_i$ for $1\le i\le N-1.$  Thus,   the smallest eigenvalue of $-L_\aa^{(N)}|_{\scr Q_1}$ is $\aa_{N+1}.$

It is easy to see that
\beg{equation*}\beg{split} &\mu^{(N)}_{\aa}(x_i):= \int_{\DD^{(N)}} x_i\mu^{(N)}_{\aa}(\d x)=\ff{\GG(\bar \aa)\GG(\aa_i+1)}{\GG(|\aa|_1+1)\GG(\aa_i)}=\ff{\aa_i}{|\aa|_1},\\
&\mu^{(N)}_{\aa}(x_i^2)= \ff{\GG(\bar \aa)\GG(\aa_i+2)}{\GG(|\aa|_1+2)\GG(\aa_i)}=\ff{\aa_i(\aa_i+1)}{|\aa|_1(|\aa|_1+1)},\ \ \  \ 1\le i\le N-1;\\
&\mu^{(N)}_{\aa}(x_ix_j)= \ff{\GG(\bar \aa)\GG(\aa_i+1)\GG(\aa_j+1)}{\GG(|\aa|_1+2)\GG(\aa_i)\GG(\aa_j)}=\ff{\aa_i\aa_j}{|\aa|_1(|\aa|_1+1)},\ \ 1\le i\ne j\le N-1.\end{split}\end{equation*}
Then
\beg{equation*}\beg{split} &\mu^{(N)}_{\aa}(u_i)= \ff 1 {|\aa|_1} \sum_{k=1}^N \theta_{ik}\aa_k=0,\ \ 1\le i\le N-1;\\
& \mu^{(N)}_{\aa,\ll}(u_N)= \sum_{i=1}^N\ff{\aa_i}{|\aa|_1} -\ff{\tt\aa}{|\aa|_1}=0.\end{split}\end{equation*} So, $\{u_i\}_{1\le i\le N}\subset \scr Q_1.$ Moreover, for $1\le i\ne j\le N-1$,
\beg{equation*}\beg{split} \mu^{(N)}_{\aa}(u_iu_j)&= \ff 1{|\aa|_1(|\aa|_1+1)}\bigg(\sum_{1\le k\le N} \theta_{ik}\theta_{jk}\aa_k(\aa_k+1) +\sum_{1\le k\ne l\le N} \theta_{ik}\theta_{jl}\aa_k\aa_l\bigg)\\
&=\ff 1 {|\aa|_1(|\aa|_1+1)}\bigg\{\Big(\sum_{1\le k\le N} \theta_{ik}\aa_k\Big)\sum_{1\le l\le N} \theta_{jl}\aa_l +\sum_{1\le k\le N} \theta_{ik}\theta_{jk}\aa_k\bigg\}=0,\end{split}\end{equation*}
and for any $1\le i\le N-1,$
\beg{equation*}\beg{split} &\mu^{(N)}_{\aa}(u_iu_N)  = \sum_{1\le k,j\le N} \theta_{ij} \mu(x_jx_k)\\
&= \ff 1 {|\aa|_1(|\aa|_1+1)} \sum_{1\le j\le N} \theta_{ij} \aa_j (\aa_j+1) + \ff 1 {|\aa|_1(|\aa|_1+1)} \sum_{1\le k\ne j\le N}\theta_{ij} \aa_j\aa_k\\
&= \ff 1 {|\aa|_1(|\aa|_1+1)} \sum_{1\le k,j\le N}\theta_{ij}\aa_j\aa_k + \ff 1 {|\aa|_1(|\aa|_1+1)} \sum_{1\le j\le N}\theta_{ij}\aa_j=0.\end{split}\end{equation*} Since   $\{\theta_i\}_{i=1}^{N-1}$ is a basis of $\R^{N-1}$, we have
$${\rm dim\ span} \{u_i: 1\le i\le n-1\}= N-1= {\rm dim}\, \scr Q_1.$$ In conclusion,   $\{u_i\}_{1\le i\le N}$ is an orthogonal basis of $\scr Q_1.$

Finally,   we have
$$L_\aa^{(N)} u_i(x)= \sum_{j=1}^N (\aa_jx_{N+1}-\aa_{N+1}x_j) \theta_{ij}=-\aa_{N+1}u_i,\ \ 1\le i\le N-1,$$ and
$$L_\aa^{(N)} u_N(x)= \sum_{j=1}^N(\aa_jx_{N+1}-\aa_{N+1} x_j) =-|\aa|_1 \sum_{j=1}^Nx_j+\sum_{j=1}^N\aa_j = -|\aa|_1 u_N(x).$$
Therefore, the proof is finished. \end{proof}

\section{The whole spectrum of $L_\aa^{(N)}$}

\def\cH{\scr H} \def\cF{\scr F} \def\cG{\scr G}\def\cP{\scr P} \def\cK{\scr K}

For $d\in\Z_+$, let $\cH_d$ be the space of homogeneous polynomials of total degree $d$ in the variables $x_1$, ..., $x_N$.
Denote by $\tt\pi_d$ the natural projection from $\PP_\infty$ to $\cH_d$ which only keeps the $d$-homogeneous part of a polynomial.
Let $L^{(N)}_{\aa,d} = (\tt\pi_d L_\aa^{(N)})|_{\cH_d}$ be the restriction of the operator $\tt\pi_d  L_\aa^{(N)}$ to $\cH_d$ and denote $-\Lambda_d$ its spectrum, seen as a multi-set (namely with multiplicities).
From the above considerations, the spectrum $\Lambda$ of $-L_\aa^{(N)}$ is equal to $\cup_{d\in \Z_+} \Lambda_d$, as a multi-set.
 We can write
$$
\bar L^{(N)}_{\aa,d} =-|\cdot|_1\tt L^{(N)}_{\aa,d} -\alpha_{N+1}\hat  L^{(N)}_{\aa,d} ,$$
where $\tt L^{(N)}_{\aa,d} :  \cH_d\to \cH_{d-1}$ and $\hat  L^{(N)}_{\aa,d} : \cH_d\to \cH_{d}$ are respectively the restriction to $\cH_d$ of the operators
$$
\tt L_\aa^{(N)} :=  \sum_{1\le n\le N} \big(x_n\pp_n^2+\alpha_n\pp_n\big),\ \ \
\hat  L_\aa^{(N)}:=  \sum_{1\le n\le N} x_n\pp_n.$$
The crucial point of the previous decomposition is that $\hat  L^{(N)}_{\aa,d}=d  I_{\cH_d}.$
Denote by $\tt \Lambda_d$ the spectrum of $|\cdot|_1\tt L^{(N)}_{\aa,d} $, we thus have
$$
\Lambda_d=\tt \Lambda_d +d\alpha_{N+1}.$$
Note that $\Lambda_0=\tt \Lambda_0=\{0\}$. The next result enables to compute by iteration
$\tt \Lambda_d$ for all $d\in\Z_+$.

\begin{prp}\label{pro2}
For any $d\in\Z_+$, we have
$$
\tt \Lambda_{d+1}=(2d+\tt\alpha +\tt \Lambda_d)\cup\{0[C(N,d+1)-C(N,d)]\},$$
where $\{0[l]\}$ is the multi-set with 0 repeated $l$ times, for $l\in\Z_+$ (more generally $[l]$ will stand for the multiplicity $l$), and where
$C(N,d)$ is the dimension of $\cH_d$, namely
$$
C(N,d)=\binom{d+N-1}{d}.$$
\end{prp}
\beg{proof}
Consider $\lambda\in \tt\Lambda_{d+1}$ and let $\varphi\in\cH_d$ be an associated eigenvector (non-zero).
We have
$$
|\cdot|_1\tt L^{(N)}_{\aa,d+1}\varphi = \lambda \varphi.$$
Since $L^{(N)}_{\aa,d+1}[\varphi]$ belongs to $\cH_{d}$, there are two possibilities:
either $\lambda=0$, or $\varphi=|\cdot|_1\psi$
for some  $\psi  \in \cH_d$ such that
\beq\label{*W}
\tt L^{(N)}_{\aa,d+1}(|\cdot|_1\psi) = \lambda \psi.\end{equation}
We consider the latter situation, since the former case leads to the multi-set $\{0[C(N,d+1)-C(N,d)]\}$.

We compute at point $x$ that
\beg{equation}\label{*W2}\beg{split}
{\tt L}^{(N)}_{\aa,d+1}(|\cdot|_1\psi)&=  |x|_1{\tt L}^{(N)}_\aa   \psi +  \psi\tt L^{(N)}_\aa  |\cdot|_1+2\sum_{1\le n\le N}x_n \pp_n \psi \\
&= |x|_1\tt L^{(N)}_{\aa,d}\psi+ \psi\sum_{1\le n\le N}\alpha_n+2\sum_{1\le n\le N}x_n\pp_n\psi\\
 &= |x|_1{\tt L}^{(N)}_{\aa,d} \psi +(\tt\alpha +2d)\psi.\end{split}\end{equation}
So, it follows from \eqref{*W} that $\lambda-\tt\alpha -2d$ is an eigenvalue of the operator $|\cdot|_1\tt L^{(N)}_{\aa,d}$,
namely belongs to $\tt \Lambda_d$. Thus, $$\tt \LL_{d+1}\subset (2d+\tt\alpha +\tt \Lambda_d)\cup\{0[C(N,d+1)-C(N,d)]\}.$$

On the other hand,  if $\ll'\in \tt \LL_d$ then $|\cdot|_1\tt L^{(N)}_{\aa,d}\psi= \ll' \psi$ for some $0\ne \psi\in \scr H_d$. Then \eqref{*W2} implies
$$\tt L^{(N)}_{\aa,d+1}(|\cdot|_1\psi)= |\cdot|_1\tt L^{(N)}_{\aa,d}\psi +(\tt\alpha +2d)\psi= (\ll'+\tt\aa+2d)\psi.$$
Therefore, $\ll'+\tt\aa+2d\in \tt\LL_{d+1}$; that is, $\tt \LL_{d+1}\supset (2d+\tt\alpha +\tt \Lambda_d).$
Then the proof is finished.  \end{proof}

The previous arguments amount to an iterative construction of the eigenvectors:
for any $d\in\Z_+$, let
$\tt\cF_d$ be the set of eigenvectors of $\tt L^{(N)}_{\aa,d} $ and $\cG_d$ be
% the set of polynomial functions from $\cH_d$ which don't admit $y_{N}$
% as a factor and which are in
 the kernel of $\tt L^{(N)}_{\aa,d} $.
Then we have
$$\forall d\in\Z_+,\qquad \tt\cF_{d+1}= \cG_{d+1}\cup y_N\tt\cF_d.$$
Indeed, in the above proof, functions $\varphi\in\tt\cF_{d+1}$ of the form $y_N\psi$ with $\psi \in\tt\cF_d$
are associated to eigenvalues of the form $\tt\alpha+2d+\lambda$, where $\lambda\in\tt\Lambda_d$.
From Lemma \ref{L1}, we know that $\lambda\geq 0$, so that $\tt\alpha+2d+\lambda>0$
and $\varphi$ does not belong to the kernel of $\tt L^{(N)}_{d+1}$.
Conversely, we have seen that all the other eigenvectors belong to the kernel of $\tt L^{(N)}_{d+1}$.
Thus we get the following characterization of the kernel of $\tt L^{(N)}_{\aa,d} $: it consists exactly into the eigenvectors
of $\tt L^{(N)}_{\aa,d} $ which  don't admit $y_{N}$
 as a factor.

Note that $\tt \cF_d$ is also the set of eigenvectors of $L^{(N)}_{\aa,d} $.
To get the eigenvectors of our initial operator $L$, we construct by iteration on $d\in\Z_+$ the following subsets $\cF_d$ of $\cP_d$.
First we take $\cF_0:=\tt\cF_0=\cP_0$.
Next, if $\cF_d$ has been constructed, then for any $f\in\tt\cF_{d+1}$, there exists a unique $g_f\in\cP_d$ such that
$f+g_f$ is orthogonal to $\cP_d$ in $\LL^2(\mu)$. Then we define
$$
\cF_{d+1}:=\{f+g_f:  f\in \tt\cF_{d+1}\}.$$
The set of eigenvectors of $L$ is $\cup_{d\in\Z_+}\cF_d$.

From Proposition \ref{pro2}, it is possible to parametrize the spectrum $\Lambda$ of $-L$ in the following way.
Let $\cK$ be the set of elements of the form $(k_1,k_2, ..., k_{r},k_{r+1})$, where $r\in\Z_+$
and $0\leq k_1<k_2<\cdots <k_r<k_{r+1}$.
Define a mapping $K:  \cK\to \Lambda$ via
$$ \forall k:= (k_1,k_2, ..., k_{r},k_{r+1})\in\cK,\qquad K(k):= 2(k_1+\cdots k_r)+ r\tt\alpha+ k_{r+1}\alpha_{N+1}.$$
Then $K$ is surjective. It is truly one-to-one, if and only if  1, $\tt \alpha$ and $\alpha_{N+1}$ are independent when $\R$
is seen as a vector space over $\Q$. Let us call this situation  generical over the choice of the parameters $\alpha:= (\alpha_n)_{1\le n\le N+1}$.

The multiplicities can also be recovered.
Consider the mapping $D: \cK\to \N$ defined by
$$ D(k):=\sum_{1\le l\le r} C(N, k_l) +\sum_{1\le l\le k_{r+1}-1, l\notin\{k_1, k_2, ..., k_r\}}\big\{C(N,l+1)-C(N,l)\big\}$$ for $
  k:= (k_1,k_2, ..., k_{r},k_{r+1})\in\cK.$
Then the multiplicity of an eigenvalue $\lambda\in \Lambda$
is given by
$$
\sum_{k\in K^{-1}(\lambda)} D(k).$$
In particularly, generically, we have $\Lambda=\{K(k)[D(k)]: k\in \cK\}$.

\section{Proofs of Theorems \ref{T1.2} and \ref{T1.3}}

To prove the first assertion, let $W$ be the $L^1$-Wasserstein distance induced by $\rr(x,y):=|x-y|_1$ on $\scr P(\DD^{(\infty)}),$ the set of all
  probability measures on $\DD^{(\infty)}$. That is, for any $\mu,\nu\in \scr P(\DD^{(\infty)}),$
  $$W(\mu,\nu):= \int_{\pi\in \scr(\mu,\nu)} \int_{\DD^{(\infty)}\times \DD^{(\infty)}} |x-y|\pi(\d x,\d y),$$ where $\scr C(\mu,\nu)$ is the set of all couplings for $\mu$ and $\nu$; i.e. $\pi\in \scr(\mu,\nu)$ if and only if it is a probability measure
  on $\DD^{(\infty)}\times \DD^{(\infty)}$
  such that $$\pi(\d x\times \DD^{(\infty)})= \mu(\d x),\ \ \pi(\DD^{(\infty)}\times\d y)= \nu(\d y).$$
It is well known that the metric $W$ is complete and induces the weak topology on $\scr P(\DD^{(\infty)}),$ see e.g. \cite[Theorems 5.4 and 5.6]{Chen}.
So, for the proof of Theorem \ref{T1.2} we only need to show that   $\{\mu_{\aa,\aa_\infty}^{(n)}\}_{n\ge 1}$ is $W$-Cauchy sequence.

\beg{proof}[Proof of Theorem \ref{T1.2}] (1) To prove that $\{\mu_{\aa,\aa_\infty}^{(n)}\}_{n\ge 1}$ is a $W$-Cauchy sequence, we use
the partition property of the Dirichlet distribution mentioned in Section 1. For any $n> m\ge 1$, let $(X_1,\cdots, X_{n+1})$ have law
$\tt\mu_{\aa^{(n)}}^{(n+1)}$. By the partition property, $\big(X_1,\cdots, X_{m-1}, \sum_{i=m}^nX_i, X_{n+1}\big)$ has law $\tt\mu_{\aa^{(m)}}^{(m+1)}$.
So,  $(X_1,\cdots,X_{m-1}, \sum_{i=m}^nX_i)$ has law $\mu_{\aa^{(m)}}^{(m)}$ while  $(X_1,\cdots,X_n)$ has law $\mu_{\aa^{(n)}}^{(n)}$. Thus,   the laws of $(X_1,\cdots, X_{m-1}, \sum_{i=m}^nX_i,0,0,\cdots,0)$ and $(X_1,\cdots, X_n,0,0,\cdots,0)$ are $\mu_{\aa,\aa_\infty}^{(m)}$ and $\mu_{\aa,\aa_\infty}^{(n)}$ respectively. Then, by the definition of $W$ and noting that $|\aa|_1<\infty$, we have
\beg{equation*}\beg{split} &\limsup_{m\to\infty}\sup_{n\ge m+1} W(\mu_{\aa,\aa_\infty}^{(m)}, \mu_{\aa,\aa_\infty}^{(n)})\le 2\limsup_{m\to\infty}\sup_{n\ge m+1}\sum_{i=m+1}^n \E|X_i|\\
 &= \limsup_{m\to\infty}\sup_{n\ge m+1}\sum_{i=m+1}^n \ff{2\aa_i}{\aa_\infty+ \sum_{i=n+1}^\infty \aa_i}=0.\end{split}\end{equation*}  Therefore, $\{\mu_{\aa,\aa_\infty}^{(n)}\}_{n\ge 1}$ is a $W$-Cauchy sequence and the proof of the first assertion is finished.

(2) It suffices to prove
\beq\label{CL} \EE^{(\infty)}_{\aa,\aa_\infty}(f,g)= - \int_{\DD^{(\infty)}} (f L_{\aa,\aa_\infty}^{(\infty)} g) \,\d\mu_{\aa,\aa_\infty}^{(\infty)},\ \ f,g\in \scr FC^2.\end{equation}
For any $f,g\in \scr FC^2$, there exist $m\in\N$ and $f_m,g_m\in C^2(\R^m)$ such that
$$f(x)= f_m(x_1,\cdots, x_m),\ \ \ g(x)= g_m(x_1,\cdots, x_m), \ \ x\in \DD^{(\infty)}.$$
So, by the definition of $\mu_{\aa,\aa_\infty}^{(n)}$ and using \eqref{E2'}, we have
\beq\label{GG}-\int_{\DD^{(\infty)}} (f  L_{\aa^{(n)}}^{(n)} g)\d \mu_{\aa,\aa_\infty}^{(n)}=   \int_{\DD^{(\infty)}}\bigg\{\Big(1-\sum_{1\le i\le n} x_i\Big) \sum_{i=1}^m  x_i(\pp_i f)(\pp_ig)\bigg\}\d\mu_{\aa,\aa_\infty}^{(n)}.\end{equation}
Since $\mu_{\aa,\aa_\infty}^{(n)}\to\mu_{\aa,\aa_\infty}^{(\infty)}$ weakly, and it is easy to see that
\beg{equation*}\beg{split}
&\lim_{n\to\infty} \sup_{x\in \DD^{(\infty)}}|f  L_{\aa^{(n)}}^{(n)} g -fL^{(\infty)}_{\aa,\aa_\infty}g|(x)=0,\\
&\lim_{n\to\infty} \sup_{x\in \DD^{(\infty)}}\bigg|\Big(1-\sum_{1\le i\le n} x_i\Big) \sum_{i=1}^m  x_i(\pp_i f)(\pp_ig)
-\Big(1-\sum_{i=1}^\infty x_i\Big) \sum_{i=1}^m  x_i(\pp_i f)(\pp_ig)\bigg|=0,\end{split}\end{equation*}
by letting $n\to\infty$ in \eqref{GG} we prove \eqref{CL}.

(3) Finally, as was shown in (2) that the desired Poincar\'e inequality follows by applying  Theorem \ref{T1.1} to $\mu_{\aa^{(n)}}^{(n)}$ on
$\DD^{(n)}$ then letting $n\to\infty.$  So, $\gap(L_{\aa,\aa_\infty}^{(\infty)})\ge \aa_\infty.$ On the other hand, let
$$u(x)= \aa_2 x_1-\aa_1 x_2,\ \ \ x\in\DD^{(\infty)}.$$ We have
$$L_{\aa,\aa_\infty}^{(\infty)}u(x)= \big\{\aa_1(1-|x|_1) -\aa_\infty x_1\big\}\aa_2 - \big\{\aa_2(1-|x|_1) -\aa_\infty x_2\big\}\aa_1
=-\aa_\infty u(x),\ \ x\in \DD^{(\infty)}.$$ This implies $\gap(L_{\aa,\aa_\infty}^{(\infty)})\le \aa_\infty.$ In conclusion, we have
$\gap(L_{\aa,\aa_\infty}^{(\infty)})= \aa_\infty.$\end{proof}

\beg{proof}[Proof of Theorem \ref{T1.3}] (a) For the first assertion, we only need to prove that $\{P_{x,T}^{(n)}\}_{n\ge 1}$ is a Cauchy sequence with respect to the $L^1$-Wasserstein distance
$$W_T(P, P'):= \inf_{\Pi\in \scr C(P,P')} \int_{\OO_T\times\OO_T} \|\xi-\eta\|_{1,\infty}  \Pi(\d \xi,\d \eta).$$
To this end, for any $n>m\ge 2$, we construct a coupling of $P_{x,T}^{(n)}$ and $P_{x,T}^{(m)}$ as follows.

Firstly, let $(X_i^{(n)}(t))_{1\le i\le n}$ solve the following SDE with $X^{(n)}_0=x^{(n)}$:
\beq\label{FY} \beg{split} \d X_i^{(n)}(t)=&\Big[\aa_i\big(1-|X^{(n)}(t)|_1\big)-\aa_\infty X_i^{(n)}(t)\Big]\d t \\ & +\ss{2(1-|X^{(n)}(t)|_1)X_i^{(n)}(t)}\, \d B_i(t),\ 1\le i\le n-1;\\
 \d X_n^{(n)}(t)=&\Big[\sum_{j=n}^\infty \aa_j \big(1-|X^{(n)}(t)|_1\big)-\aa_\infty X_n^{(n)}(t)\Big]\d t \\ & +\ss{2(1-|x^{(n)}(t)|_1)X_n^{(n)}(t)}\, \d B_n(t),\ \ t\in [0,T],\end{split}\end{equation} where $(B_i(t)_{1\le i\le n}$ are independent one-dimensional Brownian motions.  Then $P_{x,T}^{(n)}$ is the distribution of $(X^{(n)}(t))_{t\in [0,T]}.$

Next, let
\beq\label{FY2} X^{(m)}_i(t)= X^{(n)}_i(t)\ \text{for}\  1\le i\le m-1,\   \text{and\ } X_m^{(m)}(t)= \sum_{j=m^n} X_j^{(n)}(t),\ \ t\in [0,T].\end{equation}
Then $X^{(m)}(0)=x^{(m)}$ and by \eqref{FY},
\beg{equation*}  \beg{split} \d X_i^{(m)}(t)=&\Big[\aa_i\big(1-|X^{(m)}(t)|_1\big)-\aa_\infty X_i^{(m)}(t)\Big]\d t \\ &+\ss{2(1-|x^{(m)}(t)|_1)X_i^{(m)}(t)}\, \d B_i(t),\ 1\le i\le m-1;\\
 \d X_m^{(m)}(t)=&\Big[\sum_{j=m}^\infty \aa_j \big(1-|X^{(m)}(t)|_1\big)-\aa_\infty X_m^{(m)}(t)\Big]\d t \\
 &+\ss{2(1-|x^{(m)}(t)|_1)X_m^{(m)}(t)}\, \d \tt B_m(t),\ \ t\in [0,T],\end{split}\end{equation*} where
$ \d\tt B_m(t):= \ff 1 {\ss{X_m^{(m)}(t)}}\sum_{i=m}^n \ss{X_i^{(n)}(t)}\,\d B_i(t)$ is a one-dimensional Brownian motion independent of $(B_i(t))_{1\le i\le m-1}$. Therefore, $(X^{(m)}(t))_{t\in [0,T]}$ has law $P_{x,T}^{(m)}$.

Now, by \eqref{FY2} and the definition of $W_T$, we have
\beq\label{FY3} W_T(P_{x,T}^{(n)},  P_{x,T}^{(m)}) \le \E \sup_{t\in [0,T]} |X^{(m)}(t)- X^{(n)}(t)|_1= \E\sup_{t\in [0,T]}\sum_{j=m+1}^n X_j^{(n)}(t).\end{equation}
Let $Z(t)= \sum_{j=m+1}^n X_j^{(n)}(t).$ By \eqref{FY} we have
$$\d Z(t) \le \Big(\sum_{j=m+1}^\infty \aa_j\Big) \d t +\sum_{j=m+1}^n \ss{s(1-|X^{(n)}(t)|_1)X_i^{(n)}(t)}\, \d B_i(t).$$ So,
$$ Z(t)\le \sum_{j=1+m}^\infty (x_j + t \aa_j) + \sum_{j=m+1}^n \int_0^t \ss{s(1-|X^{(n)}(s)|_1)X_i^{(n)}(s)}\, \d B_i(s)=:\bar Z (t),\ \ t\in [0,T].$$
Since $Z(t)\ge 0$, $\bar Z(t)$ is a nonnegative submartingale. Then by Kolmogorov's inequality,
$$\P\Big(\sup_{t\in [0,T]} Z(t)\ge \ll\Big)\le \P\Big(\sup_{t\in [0,T]} \bar Z(t)\ge \ll\Big)\le \ff 1 \ll \E\bar Z(T) =\ff 1 \ll \sum_{j=m+1}^\infty(x_j+\aa_jT),\ \ \ll>0.$$
Since $Z(t)\le 1$, this implies
$$\E \sup_{t\in [0,T]}Z(t) \le \ll + \P\Big(\sup_{t\in [0,T]} Z(t)\ge \ll\Big)\le \ll+ \ff 1 \ll \sum_{j=m+1}^\infty(x_j+\aa_jT),\ \ \ll>0.$$
Taking $\ll= \ss{\sum_{j=m+1}^\infty(x_j+\aa_jT)},$ and combining with \eqref{FY3},
we obtain
$$\lim_{m\to\infty}\sup_{n\ge m+1} W_T(P_{x,T}^{(n)},  P_{x,T}^{(m)}) \le 2 \lim_{m\to\infty} \ss{\sum_{j=m+1}^\infty(x_j+\aa_jT)}=0.$$
Therefore, the first assertion is proved.

(b) Let $f\in \F C^2$. We have $f(x)= f(x_1,\cdots, x_m)$ for some $m\ge 1$ and $f\in C^2(\DD^{(m)})$. For the coordinate process $X(t)$,   define
$$M^{(n)}(t) =f(X(t))-\int_0^t L_{\aa,\aa_\infty}^{(n)} f(X(s))\d s,\ \ n\ge m, t\in [0,T].$$ Then $(M^{(n)}_t)_{t\in [0,T]}$ is a $P_{x,T}^{(n)}$-martingale; that is, for any $0<s<t\le T$, and any bounded Lipschitz continuous function $g $ on $\OO_T$ measurable with respect to $\F_s$,
\beq\label{WF4}\int_{\OO_T} M^{(n)}(t)(\oo) g(\oo)\d P_{x,T}^{(n)}= \int_{\OO_T} M^{(n)}(s)(\oo) g(\oo)\d P_{x,T}^{(n)}.\end{equation}
We intend to prove the same equality for $P_{x,T}^{(\infty)}$ and
$$M^{(\infty)}(t) :=f(X(t))-\int_0^t L_{\aa,\aa_\infty}^{(\infty)} f(X(s))\d s,\   t\in [0,T].$$ By an approximation argument, we may and do assume that
$f\in C_b^3(\DD^{(m)})$. In this case, $M^{(n)}(t)$ is bounded and Lipschitz on $\OO_T$ uniformly in $n\ge m$ and $t\in [0,T]$. Since $g$ is bounded and Lipschitz on $\OO_T$ as well, there exists a constant $C>0$ such that
 $$|(M^{(n)}(t)g)(\xi)-(M^{(n)}g)(t)(\eta)|\le C \|\xi-\eta\|_{1,\infty},\ \ n\ge m, \xi,\eta\in \OO_T, t\in [0,T].$$
 Therefore,
 $$\bigg|\int_{\OO_T} M^{(n)}(t)  g \d P_{x,T}^{(n)}- \int_{\OO_T} M^{(n)}(t)  g \d P_{x,T}^{(\infty)}\bigg|
 \le CW_T(P_{x,T}^{(n)},P_{x,T}^{(\infty)}),\ \ n\ge m, t\in [0,T].$$
 Combining this with \eqref{WF4}, $\lim_{n\to\infty} W_T(P_{x,T}^{(n)},P_{x,T}^{(\infty)})=0,$   $\lim_{n\to\infty} M^{(n)}= M^{(\infty)}$ and
 noting that $\{M^{(n)}g\}_{n\ge m}$ are uniformly bounded,
 we conclude that
 \beg{align*} &\bigg|\int_{\OO_T} \big[M^{(\infty)}(t) - M^{(\infty)}(s) \big]g \,\d P_{x,T}^{(\infty)}\bigg|\\
 &= \lim_{n\to\infty} \bigg|\int_{\OO_T} \big[M^{(n)}(t) - M^{(n)}(s) \big]g \,\d P_{x,T}^{(\infty)}\bigg|\\
 &\le 2C \limsup_{n\to\infty} W_T(P_{x,T}^{(n)},P_{x,T}^{(\infty)})=0.\end{align*}
 Then the proof is finished.
\end{proof}
\section{  A Discrete Model}

For any $N\geq 1$, $M\geq N+1$, consider a population of $M$ individuals of $N+1$ different types.  Divide the population into two groups: group I of types $1, \ldots, N$ and group II of type $N+1$. Focusing on group I and treat group II as outsiders or external sources.   Initially the number of type $i$ individuals is $m_i, i=1,\ldots,N+1$. The group I evolves as follows: a type $i$ individual independent of all others will wait for an exponential time at rate $\alpha_{N+1}$ and at the end of the waiting  emigrates to the outside becoming type $N+1$;  an outsider will independently wait an exponential time with rate $\alpha_i$ and immigrate to group I becoming type $i$; in addition to emigration and immigration, each couple between a type I and a type II waits for
 an exponential time with rate $2$ and when the clock rings,   either the group I individual moves out becoming an outsider or the group II individual  moves in becoming the type of the selected individual in group I.

Let $X(t)=M^{-1} (M_1(t), \ldots, M_N(t))$ denote the relative frequencies of individuals of different types in group $I$ among the whole population at time $t$. For $\aa\in (0,\infty)^{N+1}$, we  construct $X(t)$ as a multivariate Markov chain with generator
\begin{eqnarray*}
\cA_{M,\aa}^{(N)} f(x)
&=&M\sum_{i=1}^N\Big\{\alpha_{N+1}x_i \Big[f\Big(x-\frac{e_i}{M}\Big)- f(x)\Big]+\alpha_i (1-|x|_1)\Big[f\Big(x+\frac{e_i}{M}\Big)-f(x)\Big]\Big\}\\
&& +M^2\sum_{i=1}^N  (1-|x|_1)x_i\Big\{f\Big(x-\frac{e_i}{M} \Big)+f\Big(x+\frac{e_i}{M}\Big)-2f(x)\Big\},\ \ f\in C^2(\DD^{(N)})
\end{eqnarray*} for $x\in \DD_M^{(N)}:=\big\{x\in \ff 1 M \Z_+^N: \ |x|_1=\sum_{i=1}^N x_i\le 1\big\}$,
where   $e_i$ is the unit vector in the $i$th direction.
Letting $M\to\infty$ and $x\to y\in \DD^{(N)}$, one gets
$\cA_{M,\aa}^{(N)} f(x) \rightarrow L_\aa^{(N)} f(y).$

%Similarly, in the infinite-dimensional case, for $\aa\in (0,\infty)^\N$ with $|\aa|_1<\infty$ and $\aa_\infty>0$, let
% \begin{eqnarray*}
%{\cal  A}_{M,\aa,\aa_\infty}^{(\infty)} f(x)
%&=&M\sum_{i=1}^\infty\Big\{\alpha_{\infty}x_i \Big[f\Big(x-\frac{e_i}{M}\Big)- f(x)\Big]+\alpha_i (1-|x|_1)\Big[f\Big(x+\frac{e_i}{M}\Big)-f(x)\Big]\Big\}\\
%&& +M^2\sum_{i=1}^\infty  (1-|x|_1)x_i\Big\{f\Big(x-\frac{e_i}{M}\Big)+f\Big(x+\frac{e_i}{M}\Big)-2f(x)\Big\},\ \ f\in \scr FC^2
%\end{eqnarray*} for $x\in \ff 1 M \Z_+^\N$ with $|x|_1=\sum_{i=1}^\infty x_i\le 1$. When $x\to y\in \DD^{(\infty)}$ we have
%$
%\cA_{M,\aa,\aa_\infty} f(x) \rightarrow L_{\aa,\aa_\infty}^{(\infty)} f(y).$

We will see that the finite Markov chain generated by $\cA_{M,\aa}^{(N)}$ on $\DD_M^{(N)}$ is reversible with respect to the probability measure $\mu_{M,\aa}^{(N)}$:
$$\mu_{M,\aa}^{(N)}(x):= \ff{[\aa_{N+1}]_{M(1-|x|_1)}}{Z\{M(1-|x|_1)\}!} \prod_{i=1}^N \ff{[\aa_i]_{Mx_i}}{(Mx_i)!},\ \ x\in \DD_M^{(N)},$$
where $[\aa]_m:=\prod_{i=0}^{m-1}(\aa+i)$ for $\aa\ge 0$ and $m\ge 1$,   $[\aa]_0:=1$, and
$$Z:= \sum_{x\in \DD_M^{(N)} } \ff{[\aa_{N+1}]_{M(1-|x|_1)}}{\{M(1-|x|_1)\}!} \prod_{i=1}^N \ff{[\aa_i]_{Mx_i}}{(Mx_i)!}$$ is the normalization.
Moreover, for $N\ge 2$, $ \cA_{M,\aa}^{(N)}$ has the same spectral gap $\aa_{N+1}$ as  $L_\aa^{(N)}.$

\beg{thm} Let $N\ge 2$. The Markov chain generated by $\cA_{M,\aa}^{(N)}$ is irreducible and reversible with respect to $\mu_{M,\aa}^{(N)}.$ Moreover, $\cA_{M,\aa}^{(N)}$  has spectral gap $\aa_{N+1}$ in $L^2(\mu_{M,\aa}^{(N)}).$ \end{thm}

\beg{proof} (a) Denote $\gg_i=\ff{e_i}M$ for $1\le i\le N$. For any $x,y\in \DD_M^{(N)}$, let
$$q_{x,y}= \beg{cases} Mx_i\aa_{N+1} +M^2x_i(1-|x|_1), &\text{if}\ y=x-\gg_i, 1\le i\le N;\\
\aa_i M(1-|x|_1) +M^2 x_i (1-|x|_1), &\text{if} \ y=x+\gg_i, 1\le i\le N;\\
0, &\text{otherwise}.\end{cases}$$ We have
$$ \cA_{M,\aa}^{(N)} f(x) =\sum_{y\in \DD_M^{(N)}} q_{xy}\big\{f(y)-f(x)\big\},\ \ x\in \DD_M^{(N)}.$$
Since $q_{x,y}>0$ when $x,y\in \DD_M^{(N)}$ with $y=x \pm \gg_i$ for $1\le i\le N$, and $\DD_M^{(N)}$ is connected by the edges $x\to x \pm \gg_i$, we see that the Markov chain is irreducible.

Next, it is well known that
  $\cA_{M,\aa}^{(N)}$ is symmetric in
$L^2(\mu_{M,\aa}^{(N)})$ if and only if
\beq\label{Z2} \mu_{M,\aa}^{(N)}(x)q_{x,y}= \mu_{M,\aa}^{(N)}(y)q_{y,x},\ \ x,y\in \DD_M^{(N)}.\end{equation}
To verify this condition, we only need to consider the following two situations.

(a1) $y=x+\gg_i$ for some $1\le i\le N.$ In this case we have $M|x|_1\le M-1,$ and by the definition of $\mu_{M,\aa}^{(N)}$,    
$$\ff{\mu_{M,\aa}^{(N)}(y)}{\mu_{M,\aa}^{(N)}(x)}= \ff{M(1-|x|_1)(\aa_i+Mx_i)}{(\aa_{N+1}+M(1-|x|_1)-1)(Mx_i+1)} =\ff{q_{xy}}{q_{yx}}.$$

(a2) $y=x-\gg_i$ for some $1\le i\le N.$ In this case we have $Mx_i\ge 1$, and by the definition of $\mu_{M,\aa}^{(N)}$,    
$$\ff{\mu_{M,\aa}^{(N)}(y)}{\mu_{M,\aa}^{(N)}(x)}= \ff{(\aa_{N+1}+M(1-|x|_1))Mx_i}{(M(1-|x|_1)+1)(Mx_i-1+\aa_i)} =\ff{q_{xy}}{q_{yx}}.$$
In conclusion, \eqref{Z2} holds and thus, $\cA_{M,\aa}^{(N)}$ is symmetric in
$L^2(\mu_{M,\aa}^{(N)})$.

(b) For any $d\in\Z_+$, consider again $\cP_d$ the space of all polynomials (in $N$ variables) whose total degree is less than or equal to
$d$.
For any $f\in\cP_d$ and $1\le i\le N$,   $x\mapsto f(x-\gg_i)-f(x)$ and $x\mapsto f(x+\gg_i)-f(x)$ are polynomials belonging to $\cP_{d-1}$, while
  $x\mapsto f(x-\gg_i)+f(x+\gg_i)-2f(x)$ is a polynomial belonging to $\cP_{d-2}$.
From the definition of  $\cA_{M,\alpha}^{(N)}$, it follows that $\cP_d$ is preserved by $\cA_{M,\alpha}^{(N)}$.
As in Section 2, we consider for $d\in\Z_+$,
\bq
\cQ_d&:=& \{f\in \cP_d\cap L^2(\mu_{M,\alpha}^{(N)}):\  \mu_{M,\alpha}^{(N)}[fg]=0, \forall  g\in\cP_{d-1}\}\eq
(with the convention $\cQ_0=\cP_0$).
Note that for $d$ large enough, $\cQ_d=\{0\}$, nevertheless,
we still have
\bq
L^2(\mu_{M,\alpha}^{(N)})&=&\bigoplus_{d\in\Z_+}\cQ_d\eq
and the $\cQ_d$ are orthogonal.
Furthermore by symmetry of $\cA_{M,\alpha}^{(N)}$ in $L^2(\mu_{M,\alpha}^{(N)})$, each of the $\cQ_d$ is preserved by $\cA_{M,\alpha}^{(N)}$.
Thus it is sufficient to study the spectral decompositions of the restrictions of $\cA_{M,\alpha}^{(N)}$ to the $\cQ_d$.
But this is exactly the same analysis as in Section 2, because there we only used the highest monomials. Indeed, note that
for all $f\in\cQ_d$ and $1\le i\le N$, $$x\mapsto f(x-\gg_i)-f(x)+\ff{\pp_i f(x)}M,\ \
 x\mapsto f(x+\gg_i)-f(x)-\ff{\pp_i f(x)}M$$ are polynomials belonging to $\cP_{d-2}$, and
$$x\mapsto f(x-\gg_i)+f(x+\gg_i)-2f(x)-\ff{\pp_i^2 f (x)}{M^2}$$ belong to $\cP_{d-3}$, where we set $\cP_k=\{0\}$ if $k<0$.
Thus, for any polynomial $f\in \cQ_d$, the polynomials  $\cA_{M,\alpha}^{(N)}f$ and  $L_\alpha^{(N)}f$ have the same highest order term (i.e. the term of degree $d$), so that these two operators have the same spectral gap.

\end{proof}

\beg{thebibliography}{99}

\bibitem{BM73} D. Blackwell,  J. B. MacQueen,
 \emph{Ferguson's distribution via P\'olya' urn schemes,}  Ann. Statist.  1(1973), 353--355.

\bibitem{BR} J. Bakosi, J. R. Ristorcelli, \emph{A stochastic diffusion process for the Dirichlet distribution,} International J. Stoch. Anal. 2013, Article ID 842981, 7 pages.

\bibitem{Chen} M.-F. Chen, \emph{From Markov Chains to Non-Equilibrium Particle Systems,} World Scientific, 1992, Singapore.
\bibitem{ConMoi69} R. J. Connor,  J. E. Mosimann,
 \emph{Concepts of independence for proportions with a generalization of the Dirichlet distribution,}  J. Amer. Statist. Assoc. 64(1969), 194-206.
 
 \bibitem{DS09} M. D\"oring, W. Stannat, \emph{The logarithmic Sobolev inequality for the Wasserstein diffusion,}  Probab. Theory Related Fields 145(2009),  189--209.

\bibitem{QQ} C. L. Epstein, R. Mazzeo, \emph{Wright-Fisher diffusion in one dimension,}  SIAM J. Math. Anal.  42(2010),  568--608.

 \bibitem{EG93} S.N. Ethier and R.C. Griffiths,
      \emph{The transition function of a Fleming-Viot process,} Ann.  Probab.  21(1993), 1571--1590.

 \bibitem{EK81} S.N. Ethier and T.G. Kurtz,
      \emph{The infinitely-many-neutral-alleles diffusion
      model,} Adv. Appl. Probab.  13(1981), 429--452.

\bibitem{EK93} S.N. Ethier and T.G. Kurtz ,
      \emph{Fleming--Viot processes in population genetics,}
       SIAM J. Control and Optimization,  31, No. 2 (1993), 345--386.

\bibitem{Ewens04} W.J. Ewens, \emph{Mathematical Population Genetics, Vol I} Springer-Verlag, New York, 2004.

 \bibitem{FW07}  S. Feng, F.-Y. Wang,  \emph{A class of infinite-dimensional diffusion processes with connection to population genetics,} J. Appl. Probab. 44(2007), 938--949.

 \bibitem{FW14}  S. Feng,   F.-Y. Wang,   \emph{Harnack Inequality and Applications for Infinite-Dimensional GEM Processes,} to appear in Potential Analysis.
 
 \bibitem{Fer73}  T. Ferguson,  \emph{A Bayesian analysis of some nonparametric problems,} Ann. Statist. 1(1973), 209--230.

  \bibitem{FV79}
      W.H. Fleming and M. Viot, \emph{Some measure-valued Markov processes in population genetics theory,}
      Indiana Univ. Math. J., {\bf 28}(1979), 817--843.
 \bibitem{IW} N. Ikeda, S. Watanabe, \emph{Stochastic Differential Equations and Diffusion Processes,} 2nd Ed. North-Holland, Amsterdam, 1989.

\bibitem{Jac01} M. Jacobsen, \emph{Examples of multivariate diffusions: time-reversibility; a Cox-Ingersoll-Ross type process.} Department of Theoretical Statistics, Preprint 6, University of Copenhagen, 2001.

   \bibitem{P1} N. L. Johnson, \emph{An approximation to the multinomial distribution, some properties and applications,} Biometrika, 47(1960), 93--102.

   \bibitem{Mi1} L. Miclo, \emph{About projections of logarithmic Sobolev inequalities,}  Lecture Notes in Math. 1801 (J. Az\'ema, M. \'Emery, M. Ledoux, M. Yor Eds), pp. 201--221, 2003, Springer.

 \bibitem{Mi2} L. Miclo, \emph{Sur l'in\'egalit\'e de Sobolev logarithmique des op\'erateurs de Laguerre \`{a} petit param\'etre, } Lecture Notes in Math. 1801 (J. Az\'ema, M. \'Emery, M. Ledoux, M. Yor Eds), pp. 222--229, 2003, Springer.

 \bibitem{P2} J. E. Mosimann, \emph{On the compound multinomial distribution, the multivariate-distribution, and correlations among proportions,}
 Biometrika, 49(1962), 65--82.
 
  \bibitem{RS09} Max-K. von Renesse and K.T. Sturm, \emph{Entropic measure and Wasserstein diffusion,}  Ann.
      Probab. 37(2009), 1114--1191. 

  \bibitem{S} W. Stannat, \emph{On validity of the log-Sobolev
      inequality for symmetric Fleming-Viot operators,}  Ann.
      Probab. 28(2000), 667--684.

\end{thebibliography}
\end{document}